\documentclass[letterpaper, 10 pt, conference]{ieeeconf}  % Comment this line out if you need a4paper

\IEEEoverridecommandlockouts                              % This command is only needed if 
\usepackage[T1]{fontenc}   %%% ADDED BY CHARLIE - ERROR OTHERWISE %%%
\usepackage{graphics} % for pdf, bitmapped graphics files
\usepackage{mathptmx} % assumes new font selection scheme installed
\usepackage{times} % assumes new font selection scheme installed
\usepackage{amsmath} % assumes amsmath package installed
\usepackage{amssymb}  % assumes amsmath package installed
\usepackage{amsbsy}

\usepackage{mathtools}

\usepackage{dsfont}
\usepackage{bm}
\usepackage{ifthen}
\usepackage[dvipsnames]{xcolor}
\usepackage{siunitx}

\usepackage{tikz}
\usepackage{tikz-3dplot}
\usetikzlibrary{positioning}
\usetikzlibrary{shapes,arrows}
\usetikzlibrary{mindmap,trees}
\usetikzlibrary{calc}
\usetikzlibrary{patterns.meta}
\usetikzlibrary{math}
\usetikzlibrary{fit}
\usepackage{circuitikz}
\usepackage[normalem]{ulem}

\newcommand{\ud}{\, \mathrm{d}}
\newcommand{\Indicator}[1]{\mathds{1}_{#1}}

\DeclareMathOperator{\Id}{\mathbb{I}d}

\definecolor{dkgreen}{rgb}{0,0.6,0}

\newtheorem{assumption}{Assumption}[section]
\newtheorem{proposition}{Proposition}[section]
\newtheorem{definition}{Definition}[section]

\title{\LARGE \bf
Optimal Control of 1D Semilinear Heat Equations with Moment-SOS Relaxations*
}

\author{Charlie Lebarbé$^{1}$, Émilien Flayac$^{1}$, Michel Fournié$^{1}$, Didier Henrion$^{2}$, Milan Korda$^{2}$ 
\thanks{*Acknowledgments: This work has been supported by the AID (Agence de l’Innovation de D\'{e}fense) from the French Ministry of the Armed Forces (Minist\`{e}re des Arm\'{e}es).}% <-this % stops a space
\thanks{$^{1}$Charlie Lebarbé, Émilien Flayac and Michel Fournié are with the Department of Complex Systems Engineering, F\'{e}d\'{e}ration ENAC ISAE-SUPAERO ONERA
Universit\'{e} de Toulouse, France, {\tt\footnotesize \{charlie.lebarbe, emilien.flayac, michel.fournie\}@isae-supaero.fr}}%
\thanks{$^{2}$Didier Henrion and Milan Korda are with LAAS-CNRS, Universit\'e de Toulouse, France, and also with the Faculty of Electrical Engineering, Czech Technical University in Prague, Czechia, {\tt\footnotesize \{henrion, korda\}@laas.fr}}%
}

\begin{document}

\maketitle
\thispagestyle{empty}
\pagestyle{empty}

%%%%%%%%%%%%%%%%%%%%%%%%%%%%%%%%%%%%%%%%%%%%%%%%%%%%%%%%%%%%%%%%%%%%%%%%%%%%%%%%
\begin{abstract}
We use moment-SOS (Sum Of Squares) relaxations to address the optimal control problem of the 1D heat equation perturbed with a nonlinear term. We extend the current framework of moment-based optimal control of PDEs to consider a quadratic cost on the control. We develop a new method to extract a nonlinear controller from approximate moments of the solution. The control law acts on the boundary of the domain and depends on the solution over the whole domain. Our method is validated numerically and compared to a linear-quadratic controller.
%\CL{A numerical study is conducted on two application cases. First, in the linear case, a linear controller is derived from the SOS relaxation and compared to the optimal L/linear-Q/quadratic R/regulator. Second, in the semilinear case, we demonstrate the effectiveness of the method in deriving a stabilizing nonlinear feedback control.}
%\MK{I'd add a one to two sentence summary of the numerical results.}
\end{abstract}
%% SHOULD IT BE LONGER??

\begin{keywords}
    Partial Differential Equation (PDE) control, polynomial optimization, Linear Matrix Inequalities (LMI).
\end{keywords}

%%%%%%%%%%%%%%%%%%%%%%%%%%%%%%%%%%%%%%%%%%%%%%%%%%%%%%%%%%%%%%%%%%%%%%%%%%%%%%%%
\section{INTRODUCTION}

The control of nonlinear Partial Differential Equations (PDEs) presents additional challenges over the linear case, due to the inherent complexity and potential for chaotic behavior in these systems. Several contemporary problems in fluid dynamics \cite{control-navier-stokes}, as well as in fluid-structure interaction (\cite{fsi-1,num-method-fsi}),
%\didier{(There are many acroynms in the paper. Is FSI used often in our text ? I would use acronyms only if they occur often, e.g. PDE, LQR..)} \CL{(Indeed, FSI is only used thrice, including this mention)}
require advanced control techniques to successfully drive the solution to the desired target state.

A first strategy to tackle that kind of problem is to consider linear control techniques such as the Linear-Quadratic Regulator (LQR). In the case of linear PDEs, it provides the optimal control law at the continuous level of the problem (see \cite{course-LQR} for instance) which precedes any domain discretization, as well as in the finite-dimensional formulation that usually follows a PDE approximation scheme (e.g., Finite Difference/Element/Volume methods). %Note that these two approaches are not equivalent as the first one (derivation of the continuous LQR first and then discretization of the space) is optimal and the second one is generally suboptimal [INSERT REFERENCE ABOUT THAT FACT???].
One classic approach to extend the LQR to the class of nonlinear PDEs is to derive the controller from the linearized equation and apply it to the initial nonlinear problem \cite{riccati-proj}. However, if the initial state is too far from the target, LQR control laws become ineffective at stabilizing the system. We illustrate this situation later in Section \ref{sec:nonlinear_case}.
%However, the LQR is limited to initial states that are "close enough" to the target state. The linearization of the PDE introduces a domain of linear stability outside which no linear controller can drive the system towards the target.

%\textbf{[TO BE COMPLETED] Here goes the bibliography on state of the art nonlinear control techniques. These methods shall be mentioned: Pontryagin Maximum Principle (PMP) \cite{PMP-ref}, Backstepping \cite{backstepping-ref}, Control Lyapunov Function (CLF) \cite{coron-CLF}, State-dependent Riccati Equation (SDRE) \cite{SDRE-ref}, Hamilton-Jacobi-Bellman (HJB) \cite{HJB-1,HJB-2}, Model Predictive Control (MPC) [find article], Dynamical Mode Decomposition (Mode de Koopman) (DMD) \cite{DMD-ref}, IA and machine learning \cite{IA-ref,IA-ref2}. A possible hierarchy among these techniques is to first cite those that Optimize/Solve and then Discretize (SD) and, second, those that Discretize and then Optimize/Solve (DS)}

More advanced methods consist in constructing nonlinear controllers for the problem of control of nonlinear PDEs. In the literature, these methods can be divided into two main categories. The first one gathers the approaches that construct a nonlinear control law from the continuous formulation of the model and then discretize the system for numerical simulations. We shall denote these methods by "Model-Control-Discretize", see \cite[Section 17.12]{quarteroni-optim-then-discrete}. Among these, the backstepping method \cite{backstepping-ref} is particularly noteworthy for transforming the original system into a more manageable form for stabilization. It is based on control Lyapunov functions \cite{coron-CLF} that ensure the stability of the target system. However, finding an appropriate control Lyapunov function is often a difficult task because there is no systematic way to construct such functions, especially for strongly nonlinear systems.
%The process can be highly problem-specific and may require significant mathematical insight.
Similarly, Koopman-based dynamic mode decomposition approaches (e.g., \cite{DMD-ref,peitz2019koopman}) use a specific change of variables, the Koopman transform, which converts the original nonlinear PDE into a linear PDE in a new set of variables. Again, the main challenge of the method relies in finding such a Koopman transform.

The second major category comprises the "Model-Discretize-Control" methods. Notable approaches within this family include methods based on Pontryagin's maximum principle \cite{PMP-ref}. Other techniques solve the Hamilton-Jacobi-Bellman equation \cite{HJB-2} or state-dependent Riccati equations \cite{SDRE-ref} to derive a stabilizing control law. These approaches tackle high-dimensional problems and often require long computation time to converge. Recent advancements in machine learning control \cite{machine-learning} provide a feasible and distinct alternative to the methods discussed above.

%These developments offer new approaches that bypass some of the limitations associated with conventional techniques, such as the need for explicit derivation of Control Lyapunov Functions (CLFs) or solving complex, high-dimensional optimality equations.

%\textbf{Here goes the bibliography on some works from Zuazua on the semilinear heat equation.}

We consider in this paper a moment-based approach that relies on the concept of occupation measures \cite{KHL-2018}, which allows to relax the problem of control of nonlinear PDEs into a Linear Program (LP) in the space of Borel measures. The main interest of this framework is its overall convexity and its versatility, as it covers a wide variety of polynomial PDEs. The method falls into the category of "Model-Control-Discretize" approaches as it operates directly on the continuous formulation of the PDE control problem, thus avoiding the need for spatio-temporal discretization of the domain. However, this advantage comes at the expense of having to solve a sequence of convex semidefinite programming problems with Linear Matrix Inequalities (LMI) constraints of increasing size. This approach has already been applied to the 1D Burgers' equation, first to extract a stabilizing distributed control \cite{KHL-2018}, and later to derive the entropy solutions for Riemann problems \cite{moment-entropy-PDEs,moment-entropy-param-dependent}. An SOS method dual to the moment method is followed in \cite{moment-integral-variational-1} and \cite{moment-integral-variational-2} to find lower bounds on integral variational problems. Lyapunov stability certificates for linear PDEs are constructed in \cite{valmorbida_reviewer} using semidefinite programming and SOS techniques. For optimal control of linear PDEs, \cite{EDP-magron} approximates measures on infinite-dimensional spaces with finite-dimensional ones, while \cite{moment-inf-dim} leverages infinite-dimensional measures to solve a semilinear heat equation with a quadratic nonlinearity.

%For optimal control of linear PDEs, \cite{EDP-magron} proposes to approximate measures on infinite-dimensional vector spaces with measures on finite-dimensional vector spaces.
%Measures supported on infinite-dimensional vector spaces are considered in \cite{moment-inf-dim} and are used to solve numerically a semilinear heat equation.

%For optimal control of linear PDEs, \cite{EDP-magron} approximates measures on infinite-dimensional spaces with finite-dimensional ones, while \cite{moment-inf-dim} leverages infinite-dimensional measures to solve a semilinear heat equation.

\begin{figure}[ht]
\centering
\resizebox{1\columnwidth}{0.8\height}{%
%\resizebox{0.49\textwidth}{!}{%

\begin{tikzpicture}[node distance=1cm, auto]
 
% Sum shape
\node[draw,
    %fill=Rhodamine!50,
    %fill=Goldenrod!80,
    minimum width=2cm,
    minimum height=1.2cm,
    text width=3cm,
    align=center
] (sum) at (0,0) {Initial nonlinear PDE.};

% Controller
\node [draw,
    %fill=Goldenrod,
    %fill=Turquoise!50,
    minimum width=2cm,
    minimum height=1.2cm,
    text width=3cm,
    align=center,
    right=1cm of sum
]  (controller) {Weak formulation.};
 
% System H(s)
\node [draw,
    %fill=SpringGreen, 
    %fill=Turquoise!50,
    minimum width=2cm, 
    minimum height=1.2cm,
    text width=3cm,
    align=center,
    right=1cm of controller
] (system) {Moment relaxation.};
 
% Sensor block
\node [draw,
    %fill=SeaGreen,
    %fill=Turquoise!50,
    minimum width=2cm, 
    minimum height=1.2cm,
    text width=3cm,
    align=center,
    below = 1cm of system
]  (sensor) {Solve the moment relaxation.};

% Numerical Simulation block
\node [draw,
    %fill=SkyBlue!50, 
    %fill=Turquoise!50,
    minimum width=2cm, 
    minimum height=1.2cm,
    text width=3cm,
    align=center,
    left = 1cm of sensor
]  (numsim) {Extract a control law.};

% End Block
\node [draw,
    %fill=Turquoise!50,
    %fill=Rhodamine!55,
    minimum width=2cm, 
    minimum height=1.2cm,
    text width=3cm,
    align=center,
    left = 1cm of numsim
]  (endblock) {PDE numerical simulation};%{Simulate numerically the PDE.};

% Dashed rectangle around the specified blocks
\node[draw, dashed, inner sep=0.5cm, fit=(controller) (system) (sensor) (numsim), label=center:{Moment overlay}] {};

% Arrows with text label
\draw[-stealth] (sum.east) -- (controller.west);
%    node[midway,above]{$e$};
 
\draw[-stealth] (controller.east) -- (system.west);
%    node[midway,above]{$u$};
    
\draw (system.east) -- ++(1,0) node[anchor=mid](output){};
 
\draw[-stealth] (output.mid) |- (sensor.east);

\draw[-stealth] (sensor.west) -- (numsim.east);

\draw[-stealth] (numsim.west) -- (endblock.east);
 
\end{tikzpicture}

}
\caption{Moment-based approach for the control of PDEs.}% \CL{Too many colors?} \didier{(I would indeed use colors only if they are meaningful. Here I do not see any reason to use yellow resp. magenta - blue seems to correspond to the "moment overlay" but this category is already clearly emphasized by the hashed box.)}}
\label{fig:steps_extract_controller}
\end{figure}
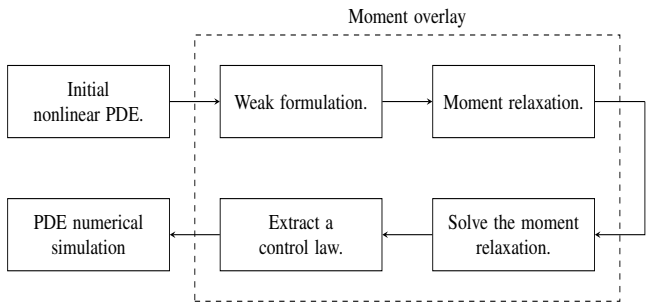

The contributions of this article are threefold: First, we extend the framework of \cite{KHL-2018} to incorporate a quadratic cost on the control, similar to the LQR. Second, we develop a new method for reconstructing a boundary control in feedback form, based on the solution over the entire domain. Third, we validate our method on a practical case where we consider the boundary control problem of the 1D heat equation perturbed with linear and nonlinear terms. Figure \ref{fig:steps_extract_controller} outlines the main steps of our method of resolution. To the best of the authors’ knowledge, the problem of control of 1D semilinear heat equations has not yet been addressed using the moments of occupation measures supported on finite-dimensional spaces. It will serve as a possible starting point for further studies on harder, higher-dimensional control problems, such as 2D semilinear heat equations, fluid-structure interaction problems, etc.

The remainder of this article is organized as follows. Section \ref{sec:section2} presents the PDE control problem that we consider and our contribution to extend the framework of \cite{KHL-2018}. Section \ref{sec:section3} discusses the validation of our method applied to a class of 1D semilinear heat equations.
%the results of our numerical simulations in the linear and nonlinear case.
We demonstrate a case where the LQR approach derived from the linearized PDE cannot control the system towards the equilibrium, whereas the moment-based approach can. Finally, Section \ref{sec:conclusion} concludes the paper and outlines future research perspectives.

%We show how to extend the framework of \cite{KHL-2018} to consider a cost function that is quadratic on the control and present a new method to extract a nonlinear control law from the moments.
%Fig. \ref{fig:steps_extract_controller} highlights the main steps of our method.

\section{Relaxation of PDE control problems into the space of occupation measures}
\label{sec:section2}

\subsection{Problem statement}
\label{sec:problem_statement}

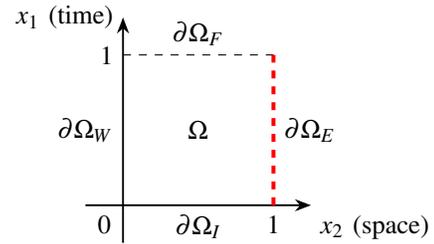
\begin{figure}[t]
\centering
%\resizebox{0.2\textwidth}{!}{%
\begin{circuitikz}
\tikzstyle{every node}=[font=\normalsize]

% Définition des variables pour les coordonnées
\def\xO{6.5}  % x de l'origine (Omega)
\def\yO{9.5}  % y de l'origine (Omega)
\def\dx{2.5}  % delta x (écart sur l'axe x)
\def\dy{2.5}  % delta y (écart sur l'axe y)

% Définir les positions des points importants
\coordinate (O) at (\xO,\yO);                 % Origine Omega
\coordinate (ystart) at (\xO, \yO - 0.5);          % Point initial sur l'axe t (temps)
\coordinate (yend) at (\xO, \yO + \dy);          % Point final sur l'axe t (temps)
\coordinate (xstart) at (\xO - 0.5, \yO);          % Point initial sur l'axe x (espace)
\coordinate (xend) at (\xO + \dx, \yO);          % Point final sur l'axe x (espace)

\coordinate (1x) at (\xO + \dx - 0.5, \yO);
\coordinate (1y) at (\xO, \yO + \dy - 0.5);

\coordinate (1x1y) at (\xO + \dx - 0.5, \yO + \dy - 0.5);

\coordinate (F1) at (\xO, \yO + 1.1*\dy);      % Point pour la frontière Omega_F
\coordinate (F2) at (\xO + \dx, \yO + \dy);    % Point pour la frontière Omega_E

% Tracé des axes
\draw [line width=0.7pt, ->, >=Stealth] (ystart) -- (yend) 
    node[pos=1,left]{$t \text{ (time)}$}
    node[pos=\dy/(\dy+0.5),left]{1}
    node[pos=(\dy + 0.5)/(2*(\dy + 0.5)),left]{$\partial \Omega_W$};

\draw [line width=0.7pt, ->, >=Stealth] (xstart) -- (xend) 
    node[pos=1,below right]{$x \text{ (space)}$} 
    node[pos=0.25/(\dx+0.5),below]{0}
    node[pos=\dx/(\dx+0.5),below]{1}
    node[pos=(\dx + 0.5)/(2*(\dx + 0.5)),below]{$\partial \Omega_I$};

% Ajouter le label Omega au centre du carré
\node at (\xO + 0.5 * \dx - 0.25, \yO + 0.5 * \dy - 0.25) {$\Omega$};

% Tracé des frontières Omega_F et Omega_E
\draw [dashed] (1y) -- (1x1y) node[pos=0.5,above]{$\partial \Omega_F$};
\draw [dashed, color=red, line width=1.5pt] (1x) -- (1x1y) node[pos=0.5,right,color=black]{$\partial \Omega_E$};

\end{circuitikz}
%}%
\caption{Definition of the domain and its boundaries. The red line indicates where the control is applied.}
\label{fig:omega_definition}
\end{figure}

We denote by $\Omega = [0,1]\!\times\![0,1]$ the time-space domain, $\partial \Omega$ its boundary and ${\scriptstyle \xi}=(t,x) \in \Omega$ the variables
%\MK{These are not state variables}
where $t$ typically represents the time and $x$ the position. We decompose the boundary of the domain into four parts $\partial \Omega = \partial \Omega_I \times \partial \Omega_F \times \partial \Omega_W \times \partial \Omega_E$. Figure \ref{fig:omega_definition} illustrates this decomposition.

We consider the problem of control of the following 1D semilinear heat equation
\begin{align}
\label{eq:heat1D}
    \begin{dcases}
        \frac{\partial y}{\partial t}({\scriptstyle \xi}) = \lambda\frac{\partial^2 y}{\partial x^2}({\scriptstyle \xi}) + \alpha y({\scriptstyle \xi}) + \eta y({\scriptstyle \xi})^3, \quad \forall {\scriptstyle \xi}\in \Omega,\\
        y({\scriptstyle \xi}) = 0, \quad \forall {\scriptstyle \xi} \in \partial \Omega_W,\\
        y({\scriptstyle \xi}) = u({\scriptstyle \xi}), \quad \forall {\scriptstyle \xi} \in \partial \Omega_E,\\
        y({\scriptstyle \xi}) = y_0({\scriptstyle \xi}), \quad \forall {\scriptstyle \xi} \in \partial \Omega_I,\\
    \end{dcases}
\end{align}
where $y$ is the unknown scalar function, $u \in L^\infty(0,1)$ is the control, $y_0 \in L^\infty(0,1)$ is a fixed initial condition, $\lambda > 0$ and $\alpha, \eta \geq 0$ are fixed parameters.
%\begin{assumption}[Finite time blow-up conditions]
%\label{assumption_blow_up}
%\hfill \break
%\emph{1) Parameters $\left(\lambda, \alpha, \eta\right)$ are such that $\alpha > \pi^2\lambda$ and $\eta > 0$.}\\
%\noindent
%\emph{2) The initial condition $y_0$ is nontrivial nonnegative.}
%\end{assumption}
%\begin{proposition}[\cite{blow_up_book}]
%\label{proposition_blow_up}
%    \emph{Under Assumption \ref{assumption_blow_up} and null boundary control, the solutions of \eqref{eq:heat1D} blow up in finite time.}
%\end{proposition}
In the next proposition, we use a weaker formulation of \cite[Proposition 3]{well_posed_pde} to state the well-posedness of system \eqref{eq:heat1D}.
\begin{proposition}
\label{proposition_well_posed}
    \emph{For all $u \in L^\infty(0,1)$, there exists a time horizon $T_u > 0$ such that system \eqref{eq:heat1D} admits a unique weak solution $y \in V_u = L^\infty((0,T_u)\times(0,1))\cap C^0([0,T_u];H^{-1}(0,1))$, where $H^{-1}(0,1)$ is the dual of $H_0^1(0,1)$.}
\end{proposition}
%\didier{The notation $L^\infty((0,T_u)\times(0,1))$ indicates that $y$ is scalar and real valued. But then it is not consistent with the notation $C^0([0,T_u];H^{-1}(0,1))$ which says that $y$ has values in $H^{-1}(0,1)$.}
%\didier{Maybe we should briefly recall what is $H^{-1}$.}\\[0.5em]
%\noindent
%\CL{\textbf{My understanding:} $H^{-1}(0,1)$ is the dual space of $H_0^1(0,1)$. We have the inclusions $H_0^1(0,1)\subset L^2(0,1) \subset H^{-1}(0,1)$ where we identify $L^2(0,1)$ with its dual, but we cannot at the same time identify $H_0^1(0,1)$ with $H^{-1}(0,1)$ (cf. Brezis Remark 3 in Chapter 5).

%My first idea would be to identify $H_0^1(0,1)$ and $H^{-1}(0,1)$ anyway with the Riesz-Fréchet representation theorem with the scalar product in $H^1$, as we do not identify $L^2(0,1)$ and its dual anywhere in the paper. However, wouldn't it mean that $y \in C^0([0,T_u];H_0^1(0,1))$ and thus $y(t,1) = 0$, which is not true? 

%\CL{My second idea would be to use the result from Brezis p.291. If $F\in H^{-1}(0,1)$, we know that there exists $f_1\in L^2(0,1)$ such that $$\forall v\in H_0^1(0,1), \quad F(v) = \int_0^1 f_1\frac{\partial v}{\partial x} \ud x.$$ Maybe we can say that we identify $F$ with $f_1$ (or its equivalence class) ????, such that $V_u$ now writes $$V_u = L^\infty((0,T_u)\times(0,1))\cap C^0([0,T_u];L^2(0,1))$$}
\begin{assumption}
\label{assumption_exist_control}
    \emph{There exists a controller $u \in L^\infty(0,1)$ such that $T_u=1$ and $y \in V_u \cap H^1(\Omega)$.}
\end{assumption}
The assumption \ref{assumption_exist_control} holds throughout this paper. The optimal control problem writes
\begin{equation}
\label{eq:optim_control_prblm}
\begin{array}{ll}
\underset{\scriptstyle \text{$(y,u) \in H^1(\Omega) \times L^\infty(0,1)$}}{\inf} & \mathcal{L}(y,u),\\
\quad \quad \quad \,\,\text{s.t.}& \text{$(y,u)$ subject to \eqref{eq:heat1D}},
\end{array}
\end{equation}
where $\mathcal{L}$ is a quadratic cost function of the form
\begin{equation}
\label{eq:cost_function}
    \mathcal{L}(y,u) = \frac{1}{2}\int_{\Omega} y({\scriptstyle \xi})^2 \ud {\scriptstyle \xi} + \frac{R}{2} \int_{0}^1 u(t)^2 \ud t
\end{equation}
and $R>0$ is fixed. We denote by $\ell_{\inf}$ the infimum of \eqref{eq:optim_control_prblm}.

\subsection{Occupation measures}
We apply the nonlinear optimal control techniques developed in \cite{KHL-2018} to solve this problem. Note that the nonlinearity we study in this article is cubic; however, our method can be immediately generalized to any other polynomial nonlinearity.
We consider $\textbf{Y} = \mathbb{R}$ such that $y$ lives in $\textbf{Y}$, $\textbf{Z} = \mathbb{R}^2$ such that $z = \nabla y$ lives in $\textbf{Z}$ and $\textbf{U} = [u_{\mathrm{min}},u_{\mathrm{max}}]$ such that $u$ lives in $\textbf{U}$.
%We consider three compact sets $\textbf{Y}$, $\textbf{Z}$ and $\textbf{U}$ such that $y$ lives in $\textbf{Y}$, $z = \nabla y$ lives in $\textbf{Z}$ and $u$ lives in $\textbf{U}$. Such sets exist because $\Omega$ is a compact set itself \CL{and $R > 0$.} \CL{It's not obvious how these sets can be determined a priori. For $\mathbf{Y}$ and $\mathbf{U}$ we should be able to get estimates but not so easy for $\mathbf{Z}$?}.
We denote by $\mu$ the occupation measure associated to $y$ on $\Omega$ and $\mu_I$, $\mu_F$, $\mu_W$ and $\mu_E$ the boundary occupation measures associated to $y$ on $\partial \Omega_I$, $\partial \Omega_F$, $\partial \Omega_W$ and $\partial \Omega_E$, respectively.
%Because we consider a quadratic cost on the control, we must adapt the framework of \cite{KHL-2018} to define $\mu_E$.
We modify the definition of the occupation measures given in \cite{KHL-2018} to include $u$ as a variable of $\mu_E$. This allows one to account for a quadratic cost on the control similar to the LQR approach.
\begin{definition}
\label{definition:occupation_measure}
\emph{The occupation measures are defined for all Borel sets $A \subset \Omega$, $A_i \subset \partial \Omega_i$ ($i \in \{I,F,W,E\}$), $B\subset \textbf{Y}$, $C\subset\textbf{Z}$ and $D\subset\textbf{U}$ by}
\begin{equation*}
\begin{split}
    &\mu(A \!\times\! B \!\times\! C) = \int_{\Omega} \Indicator{A \times B \times C}({\scriptstyle \xi},y({\scriptstyle \xi}),z({\scriptstyle \xi})) \ud {\scriptstyle \xi},\\
    &\mu_i(A_i \!\times\! B \!\times\! C) \!=\!\!\! \int_{\partial \Omega_i} \!\!\!\!\!\!\! \Indicator{A_i \times B \times C}({\scriptstyle \xi},y({\scriptstyle \xi}),z({\scriptstyle \xi})) \ud \sigma({\scriptstyle \xi}), \ i \in \{I,F,W\},\\
    &\mu_E(A_E \!\times\! B \!\times\! C \!\times\! D) \!=\!\!\! \int_{\partial \Omega_E} \!\!\!\!\!\!\!\!\Indicator{A_E \times B \times C \times D}({\scriptstyle \xi},y({\scriptstyle \xi}),z({\scriptstyle \xi}),u({\scriptstyle \xi})) \ud \sigma({\scriptstyle \xi}),
\end{split}
\end{equation*}
\emph{where $\sigma$ is the surface measure on $\partial \Omega$.}
\end{definition}
We define $\Gamma = \Omega \times \textbf{Y} \times \textbf{Z}$, $\partial \Gamma_E = \partial \Omega_E \times \textbf{Y} \times \textbf{Z} \times \textbf{U}$ and for all $i \in \{I,F,W\}$, $\partial \Gamma_i = \partial \Omega_i \times \textbf{Y} \times \textbf{Z}$. The following proposition is an immediate consequence of Definition \ref{definition:occupation_measure}.

\begin{proposition}
\label{proposition:occupation_measure}
\emph{For any bounded Borel measurable functions $h \, : \, \Gamma \rightarrow \mathbb{R}$ and $h_i \, : \, \partial \Gamma_i \rightarrow \mathbb{R}$ ($i \in \{I,F,W,E\}$), we have}
\begin{equation*}
\begin{split}
    &\int_\Omega h({\scriptstyle \xi},y({\scriptstyle \xi}),z({\scriptstyle \xi})) \ud {\scriptstyle \xi} = \int_{\Gamma} h({\scriptstyle \xi},y,z) \ud \mu({\scriptstyle \xi},y,z),\\
    &\int_{\partial \Omega_i}\!\!\!\! h_i({\scriptstyle \xi},y({\scriptstyle \xi}),z({\scriptstyle \xi})) \ud {\scriptstyle \xi} \!=\!\! \int_{\partial \Gamma_i}\!\!\! h_i({\scriptstyle \xi},y,z) \ud \mu_i({\scriptstyle \xi},y,z), i \in \{I,F,W\},\\
    &\int_{\partial \Omega_E}\!\!\! h_E({\scriptstyle \xi},y({\scriptstyle \xi}),z({\scriptstyle \xi}),u({\scriptstyle \xi})) \ud {\scriptstyle \xi} \!=\!\! \int_{\partial \Gamma_E}\!\!\! h_E({\scriptstyle \xi},y,z,u) \ud \mu_E({\scriptstyle \xi},y,z,u).
\end{split}
\end{equation*}
\emph{One should note that on the right-hand side, $y$, $z$ and $u$ are no longer functions of ${\scriptstyle \xi}$ but rather integrated variables.}
\end{proposition}
We can deduce from Proposition \ref{proposition:occupation_measure} the weak formulation of \eqref{eq:heat1D} in terms of the occupation measures (see Theorem~1 in \cite{KHL-2018} for a detailed derivation). For all test functions $\phi \in C^{\infty}(\Omega \times \textbf{Y})$, we have
\begin{subequations}
\label{eq:weak_formulation_heat1D}
    \begin{align}
    \begin{split}
    &\int_{\Gamma} \frac{\partial \phi}{\partial t} + z_1 \frac{\partial \phi}{\partial y} \ud \mu + \int_{\partial \Gamma_I} \phi \ud \mu_I - \int_{\partial \Gamma_F} \phi \ud \mu_F = 0,
    \end{split}\label{eq:weak_stokes_1} \\
    \begin{split}
    &\int_{\Gamma} \frac{\partial \phi}{\partial x} + z_2 \frac{\partial \phi}{\partial y} \ud \mu + \int_{\partial \Gamma_W} \phi \ud \mu_W - \int_{\partial \Gamma_E} \phi \ud \mu_E = 0,
    \end{split}\label{eq:weak_stokes_2} \\
    \begin{split}
    &\int_{\Gamma} \phi \left[z_1 - \alpha y - \eta y^3 \right] \ud \mu +\int_{\Gamma} \lambda\left[\frac{\partial \phi}{\partial x} + z_2 \frac{\partial \phi}{\partial y}\right] z_2\ud \mu\\
    &\quad + \int_{\partial \Gamma_W} \lambda  \phi z_2\ud\mu_W - \int_{\partial \Gamma_E} \lambda  \phi z_2\ud\mu_E = 0,
    \end{split}\label{eq:weak_domain} \\
    \begin{split}
    &\int_{\partial \Gamma_I} \phi \left[y-y_0({\scriptstyle \xi}) \right] \ud\mu_I + \int_{\partial \Gamma_E} \phi \left[y-u \right] \ud\mu_E\\
    &\quad + \int_{\partial \Gamma_W} \phi y \ud \mu_W = 0,
    \end{split}\label{eq:weak_boundary}
    \end{align}
\end{subequations}
where $\displaystyle z = [z_1 \ z_2]^T$. Equations \eqref{eq:weak_stokes_1} and \eqref{eq:weak_stokes_2} encode Stokes' formula in time and space, respectively. Equation \eqref{eq:weak_domain} is obtained by testing the PDE with test functions \(\phi\) and integrating by parts in space, while equation \eqref{eq:weak_boundary} results from testing the boundary conditions. We omitted the dependency of the occupation measures with the variables of integration ${\scriptstyle \xi}$, $y$, $z$ and $u$ in $\eqref{eq:weak_formulation_heat1D}$ for conciseness. The weak formulation allows one to write an infinite-dimensional LP whose optimal value provides a lower bound on the optimal value $\ell_{\inf}$ of \eqref{eq:optim_control_prblm}. We denote by $\mathcal{M}(A)_+$ the set of all nonnegative Borel measures with supports included in the set $A$ and by $\mathcal{K}$ the convex cone $\displaystyle \mathcal{K} = \mathcal{M}(\Gamma)_{+} \times \mathcal{M}(\partial \Gamma_I)_{+} \times \mathcal{M}(\partial \Gamma_F)_{+} \times \mathcal{M}(\partial \Gamma_W)_{+} \times \mathcal{M}(\partial \Gamma_E)_{+}$. The infinite-dimensional LP thus writes
\begin{equation}
\label{eq:LP_heat1D}
\begin{array}{ll}
\!\!\!\!\!\!\!\underset{\left(\mu, \mu_{I}, \ldots, \mu_{E}\right) \in \mathcal{K}}{\inf} \!\!\!\!\! & \displaystyle \frac{1}{2} \int_{\Gamma} y^2 \ud\mu + \frac{R}{2}\int_{\partial \Gamma_E} u^2 \ud\mu_E \\
\text { subject to } & \eqref{eq:weak_stokes_1}, \eqref{eq:weak_stokes_2}, \eqref{eq:weak_domain}, \eqref{eq:weak_boundary}, \ \forall \phi \in C^{\infty}\!(\Omega \!\times\! \mathbf{Y}), \\
& \displaystyle \int_{\partial \Gamma_i} \psi({\scriptstyle \xi}) \ud \mu_i = \int_{\partial \Omega_i} \psi({\scriptstyle \xi}) \ud\sigma({\scriptstyle \xi}),\\
&\forall \psi \in C^{\infty}(\Omega), \ \forall i \in \{I,F,W,E\}.
\end{array}
\end{equation}
We denote by $\hat{\ell}_{\inf}$ the infimum of \eqref{eq:LP_heat1D}, which satisfies $\hat{\ell}_{\inf} \leq \ell_{\inf}$. Note that the LP formulation is a relaxation of the initial problem in the sense that the set of all measures $(\mu,\mu_I,\ldots,\mu_E)$ that satisfy the previous weak formulation may be strictly larger than the set of all occupation measures corresponding to the solution(s) of the PDE. The study of the existence of a relaxation gap for the general problem of optimal control of PDEs is still an open question with the most recent results in the nonconvex case being~\cite{korda2022gap} and the convex case~\cite{henrion2024occupation,fantuzzi2022sharpness}. In our case, the existence of a relaxation gap may only have a minor influence on our results, as our goal is not to solve \eqref{eq:optim_control_prblm} exactly, but rather to present a method for constructing a stabilizing nonlinear control law in a feedback form. %\MK{If needed, the relaxation gap discussion can be removed to save space. If kept, I would suggest including some references.}

%The absence of relaxation gap was proven in \cite{no-relax-gap} for the problem of optimal control of ordinary differential equations 

\subsection{LMI relaxations}

In order to approximate the infimum $\hat{\ell}_{\inf}$ of \eqref{eq:LP_heat1D}, a hierarchy of finite-dimensional LMI relaxations is derived from the LP. The relaxation considers  polynomial test functions of the form $\displaystyle \phi(t,x,y) = t^{\alpha_1} x^{\alpha_2} y^{\alpha_3}$ in the weak formulation and truncates up to a certain relaxation degree $\displaystyle d \geq \alpha_1 + \alpha_2 + \alpha_3$. We consider $\textbf{K} \in \{\Gamma, \partial \Gamma_I, \partial\Gamma_F, \partial\Gamma_W, \partial\Gamma_E\}$ and we denote by $n_{\textbf{K}}$ the appropriate dimension such that $\textbf{K} \subset \mathbb{R}^{n_{\textbf{K}}}$. We construct a finite-dimensional outer approximation $\mathcal{M}_d(\textbf{K})_+$ of $\mathcal{M}(\textbf{K})_+$ that writes
\begin{equation}
\begin{split}
    \!\!\!\!\mathcal{M}_d(\textbf{K})_+ = \Bigl \{ \textbf{s} \in \mathbb{R}^{s_{n_\textbf{K}}(d)} \, \Big | \, &0 \preceq M_d(\textbf{s}),\\
    &0 \preceq M_d(g_i^{(\textbf{K})},\textbf{s}), \, 1 \leq i \leq m_{\textbf{K}} \Bigr \},
\end{split}
\end{equation}
where $s_{n_\textbf{K}}(d) = \binom{n_\textbf{K}+d}{d}$ is the number of monomials of $n_\textbf{K}$ variables with degree less than or equal to $d$, $\preceq$ denotes positive semidefiniteness, $\{ g_i^{(\textbf{K})}, 1\leq i \leq m_{\textbf{K}} \}$ are the polynomials that appear in the polynomial inequalities defining the semialgebraic set $\textbf{K}$, $M_d(\textbf{s})$ is the moment matrix of order $d$ and $M_d(g_i^{(\textbf{K})},\textbf{s})$ are the localizing matrices of order $d$ (for more details, see \cite{book-lasserre}). The convex cone $\mathcal{K}$ is therefore approximated by the convex semidefinite representable cone $\mathcal{K}_d = \mathcal{M}_d(\Gamma)_{+} \times \mathcal{M}_d(\partial \Gamma_I)_{+} \times \mathcal{M}_d(\partial \Gamma_F)_{+} \times \mathcal{M}_d(\partial \Gamma_W)_{+} \times \mathcal{M}_d(\partial \Gamma_E)_{+}$. The constraints of \eqref{eq:LP_heat1D} can be rewritten as a linear equation of the form $A_d\textbf{s} = b_d$ and the objective functional as a scalar product $c_d^T \textbf{s}$ where $\displaystyle \textbf{s} = \left[\textbf{s}_{\mu}^T \enspace \textbf{s}_{\mu_I}^T \enspace \cdots \enspace \textbf{s}_{\mu_E}^T \right]^T$ is the truncated vector of moments of the occupation measures up to the degree $d$. The entries of $A_d$, $b_d$ and $c_d$ depend on the coefficients of the different polynomial expressions in the weak formulation. The \textsc{Matlab} toolbox GloptiPoly 3 \cite{glopti-article} provides a well-suited framework that assembles automatically these entries. The degree $d$ finite-dimensional LMI relaxation of $\eqref{eq:LP_heat1D}$ finally writes
\begin{equation}
\label{eq:SDP_relaxation}
\begin{array}{ll}
\underset{\textbf{s} \in \mathcal{K}_d}{\inf} &c_d^T \textbf{s} \\
\,\,\,\text{s.t.} &A_d \textbf{s} = b_d.
\end{array}
\end{equation}
We solve \eqref{eq:SDP_relaxation} with the MOSEK solver based on interior-point methods \cite{sra_interior-point_2011}. It provides a sequence $\displaystyle \left(\hat{\ell}_{\inf,d}^{\,\,\mathrm{SDP}}\right)_d$ of lower bounds of $\hat{\ell}_{\inf}$ that does not decrease with the relaxation degree $d$. From the LMI relaxation, we also obtain pseudo-moments, which are approximations of the moments of the occupation measures.

\subsection{Extraction of a control law from the pseudo-moments}

After solving the LMI problem, all the pseudo-moments up to the relaxation degree $d$ have been computed. %\footnote{In fact, pseudo-moments of higher degrees may have been computed depending on the expressions of the integrands in \eqref{eq:weak_formulation_heat1D}. For instance, if we denote by $d_0$ the degree of the initial condition $y_0$, all the pseudo-moments of the measure $\mu_I$ up to degree $d+d_0$ are computed in \eqref{eq:weak_boundary}.}.
We develop in this section a method to extract from the pseudo-moments a nonlinear feedback control law that depends on the solution over the whole domain. We suppose the following form of the control
\begin{equation}
\label{eq:control_nonlinear_integral}
    u(t) = \int_0^1 \gamma(t,x,y(t,x))\ud x, \quad \forall t \in [0,1],
\end{equation}
where $\gamma$ is a multivariate polynomial of degree $m \in \mathbb{N}$. We denote by $\displaystyle \beta_{m}(t,x,y) = \left[1 \enspace t \enspace x \enspace y \enspace t^2 \enspace t x \enspace \cdots \enspace y^{m} \right]^T \in \mathbb{R}^{s_3(m)}$ the vector of all monomials up to degree $m$, sorted with the graded lexicographic ordering. There exists a vector of coefficients $c_\gamma \in \mathbb{R}^{s_3(m)}$ such that
\begin{equation}
\label{eq:gamma_scalar_product}
    \gamma(t,x,y) = \beta_m(t,x,y)^Tc_\gamma.
\end{equation}
We can multiply \eqref{eq:control_nonlinear_integral} by $\phi \in C^\infty([0,1])$, integrate in time and express the result in terms of the occupation measures. We obtain
\begin{equation}
    \int_{\partial \Gamma_E} \phi(t)u \ud \mu_E = \left( \int_{\Gamma} \phi(t) \beta_{m}(t,x,y)^T \ud \mu \right) c_\gamma,
\end{equation}
where the integral on the right-hand side is taken component-wise. The test functions are chosen as $\phi(t)=t^k$ for all $k \in \{0,\ldots,p\}$ and $p \in \mathbb{N}$, which leads to the resolution of a rectangular linear system
\begin{equation}
\label{eq:rectangular_sys}
\Phi = Bc_\gamma,
\end{equation}
where $\Phi \in \mathbb{R}^{p+1}$ and $B \in \mathbb{R}^{(p+1)\times s_3(m)}$. We have for all $k \in \{1,\ldots,p+1\}$
\begin{equation}
\label{eq:vector_phi_matrix_B}
    \Phi_k = \int_{\partial \Gamma_E} t^{k-1} u \ud \mu_E \,\text{ and }\, B_{k,\bullet} = \int_{\Gamma} t^{k-1} \beta_{m}(t,x,y)^T \ud \mu
\end{equation}
%\begin{subequations}
%    \begin{align}
%        \Phi_k &= \int_{\partial \Gamma_E} t^{k-1} u \ud \mu_E,\label{eq:vector_phi}\\
%        B_{k,\bullet} &= \int_{\Gamma} t^{k-1} \beta_{m}(t,x,y)^T \ud \mu, \label{eq:matrix_B}
%    \end{align}
%\end{subequations}
where $B_{k,\bullet}$ denotes the $k^{\mathrm{th}}$ row of the matrix $B$. Because the pseudo-moments have only been computed up to the relaxation degree $d$, \eqref{eq:vector_phi_matrix_B} imposes the condition $p\leq \min(d-1,d-m)$. 
If $p+1 > s_3(m)$, the solution $c_\gamma^*$ to the system \eqref{eq:rectangular_sys} is chosen to be the least-squares solution. Otherwise, $c_\gamma^*$ is chosen to be the minimum-norm solution for the Euclidean norm $\lVert \cdot \rVert_2$ in $\mathbb{R}^{s_3(m)}$. The choice of the norm is not inconsequential as it influences the coefficients of $\gamma$. %\MK{\sout{Without \textit{a priori} information on the form of the optimal control, other standard norms on $\mathbb{R}^{s_3(m)}$ may lead to more optimal control laws.}}
We restrict our study to the Euclidean norm and leave the investigation of other norms for future work.

Note that \eqref{eq:control_nonlinear_integral} is not the only form of the control that can be recovered from the pseudo-moments. For instance, one could consider a control that is linear in the solution
\begin{equation}
\label{eq:control_linear_integral}
    u(t) = \int_0^1 y(t,x) \gamma_{\,\ell}(t,x)\ud x,
\end{equation}
or semilinear as in
\begin{equation}
\label{eq:control_semilinear_integral}
    u(t) = \int_0^1 y(t,x) \gamma_{\,s\ell}(t,x) + y(t,x)^r \delta_{\,s\ell}(t,x)\ud x,
\end{equation}
where $t \in [0,1]$, $r \geq 2$, $\gamma_{\,\ell}$ and $\gamma_{\,s\ell}$ are multivariate polynomials of degree $m$ and $\delta_{\,s\ell}$ is a multivariate polynomial of degree $m_r$. Because they are particular cases of \eqref{eq:control_nonlinear_integral}, these control laws ultimately result in solving a rectangular linear system similar to \eqref{eq:rectangular_sys}. We will see in Section \ref{sec:linear_case} that one can be interested in a controller of the form \eqref{eq:control_linear_integral} to mimic the behavior of the LQR, while \eqref{eq:control_semilinear_integral} is interpreted in Section \ref{sec:nonlinear_case} as a perturbation of \eqref{eq:control_linear_integral} that is more robust to the cubic nonlinearity of the PDE. The key advantage of this reconstruction method is that it allows a wide range of controllers to be derived from the pseudo-moments, without the need for advanced numerical techniques beyond the resolution of \eqref{eq:rectangular_sys}.
%where $t \in [0,1]$, $r \in \mathbb{N}$ and $\delta$ is a multivariate polynomial of degree $m_r \CL{\in \mathbb{N}}$. All these control laws lead to the resolution of rectangular linear systems similar to \eqref{eq:rectangular_sys}. The advantage of this method of reconstruction lies in the fact that a wide range of controllers can be accessed from the pseudo-moments, which does not require any advanced numerical techniques beyond solving a rectangular linear system.

\section{Numerical simulations}
\label{sec:section3}

For the next two sections \ref{sec:LQR} and \ref{sec:linear_case}, we will assume that $\eta = 0$ in \eqref{eq:heat1D}. Notice that this case is the linearized version around the equilibrium state $y=0$ of \eqref{eq:heat1D} when $\eta \neq 0$. The resulting PDE is simply the heat equation, shifted by the linear term $\alpha y$. Note that the eigenvalues $\theta_k$ of the operator $\lambda \frac{\mathrm{d}^2}{\mathrm{d}x^2} + \alpha \Id$ ($\Id$ denotes the identity operator) with homogeneous Dirichlet boundary conditions are given by
\begin{equation}
    \theta_k = \alpha - \lambda \pi^2 k^2, \quad \forall k \geq 1.
\end{equation}
Thus, if $\alpha > \lambda \pi^2$, there is at least one positive eigenvalue and the solution of \eqref{eq:heat1D} naturally diverges if no control is applied.

We will first consider the LQR controller, which gives the theoretical optimal solution to the problem of Section \ref{sec:problem_statement} in the linear case, and compare it to linear and nonlinear controllers derived from the moment-based method.

\subsection{The linear-quadratic regulator}
\label{sec:LQR}
\begin{figure}[h]
    \centering
    \includegraphics[width=1\columnwidth]{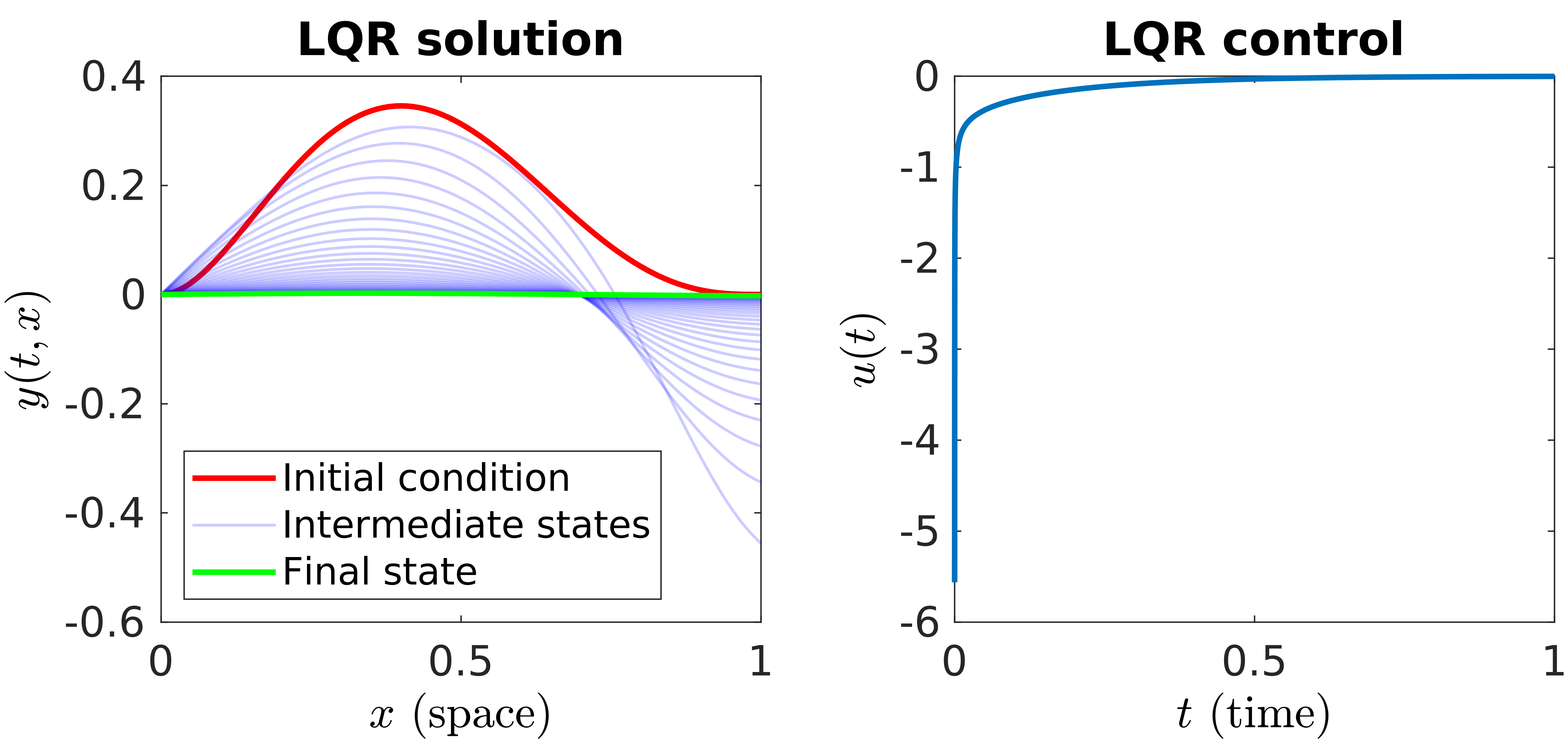}
    \caption{Numerical simulation of \eqref{eq:heat1D} with $\eta = 0$ and the LQR controller.}
    \label{fig:numsim_LQR_linear_heat}
\end{figure}

We use $P_1$-Lagrange finite elements on a uniform grid in space to approximate the solution of the PDE, which leads to the resolution of
\begin{align}
\label{eq:ODE_control}
    \begin{dcases}
        M\dot{Y}(t) = AY(t) + BU(t), \quad \forall t>0,\\
        Y(0) = Y_0 \in \mathbb{R}^N,
    \end{dcases}
\end{align}
where $Y \in \mathbb{R}^N$ approximates the solution $y(t,\cdot)$ with $N$ degrees of freedom, $U\in \mathbb{R}$ approximates the control $u(t)$ at $x = 1$, and $M,A \in \mathbb{R}^{N \times N}$ and $B \in \mathbb{R}^{N}$ are the finite elements matrices and vector. The discrete cost function eventually reads
\begin{equation}
\label{eq:discrete_cost_function}
    %\!\!\!\!\!L(Y,U) = \frac{1}{2}\int_0^{\infty} \left( Y(t)^T Q Y(t) + R U(t)^2 \right)\ud t,
    L(Y,U) = \frac{1}{2}\int_0^{\infty} Y(t)^T Q Y(t) \ud t + \frac{R}{2}\int_0^{\infty} U(t)^2 \ud t,
\end{equation}
where $R>0$ is the same as in \eqref{eq:cost_function}
and $Q \in \mathbb{R}^{N \times N}$ is a positive semidefinite weight matrix on $Y$. The optimal pair minimizing $L$ is given by $\left(\Bar{Y},\Bar{U}\right) = \left(Y,-R^{-1}B^TPMY\right)$ where $P \in \mathbb{R}^{N \times N}$ is symmetric and solution to an Algebraic Riccati Equation (ARE) (see \cite{ARE-book}).
%The resolution of the ARE requires $N\left(N+1\right)/2$ unknown coefficients to be computed, hence the difficulties encountered when considering PDEs that require very fine meshing of the domain. In practice, the LQR is not derived on the full system but rather on a projection of the Riccati equation onto the unstable eigenspaces of the system (cf. \cite{riccati-proj}). Thus, for a small number $N_u$ of unstable components, the computational complexity can be significantly reduced. However, for large systems, the derivation of the spectrum is a numerical challenge on its own. Furthermore, the controller derived from this projection is subobtimal, as the optimal controller is given by the LQR derived on the full system.
Notice that unlike in \eqref{eq:cost_function}, infinite-time integrals are considered in the discrete cost \eqref{eq:discrete_cost_function}. This choice is justified by the results of the numerical simulations, as the infinite-time horizon LQR stabilizes the solution in short time (see Fig. \ref{fig:numsim_LQR_linear_heat}).%We could be more precise and solve a differential Riccati equation on the given finite-time horizon.

We present in Fig. \ref{fig:numsim_LQR_linear_heat} the results of the numerical simulations for the LQR. The parameters are fixed to $R = 10^{-3}$, $\lambda = 0.5$ and $\alpha = 0.2 + \lambda\pi^2$ such that exactly one unstable eigenvalue is forced. We choose a polynomial initial condition $\displaystyle y_0(x) = 10\,x^2\left(1-x\right)^3$ that satisfies the homogeneous Dirichlet boundary conditions and has a nonzero projection onto the first unstable mode of the PDE. The solution of the PDE is computed using a backward differentiation scheme of order 2 with a time step $\Delta t = 10^{-4}$. The space discretization uses $P_1$-Lagrange finite elements with a uniform mesh size $h = 0.01$.

Figure \ref{fig:numsim_LQR_linear_heat} shows that the LQR successfully controls the solution to $y=0$.
%The variations of the controller are very steep around $t = 0$ because the parameter $R$ is close to 0 and allows the norm of the control in \eqref{eq:discrete_cost_function} to be large. It acts as a regularization parameter \CL{when it is increased, yielding a slower and more passive controller.}
The optimal value $L^*$ of the discrete cost function is given by the LQR solution which, in this case, yields $L^* = L_{\mathrm{LQR}} \approx \num{1.829e-3}$. In the following, every numerical result will be rounded to three decimal places.

\subsection{The linear case}
\label{sec:linear_case}

\begin{figure}[h]
    \centering
    \includegraphics[width=1\columnwidth]{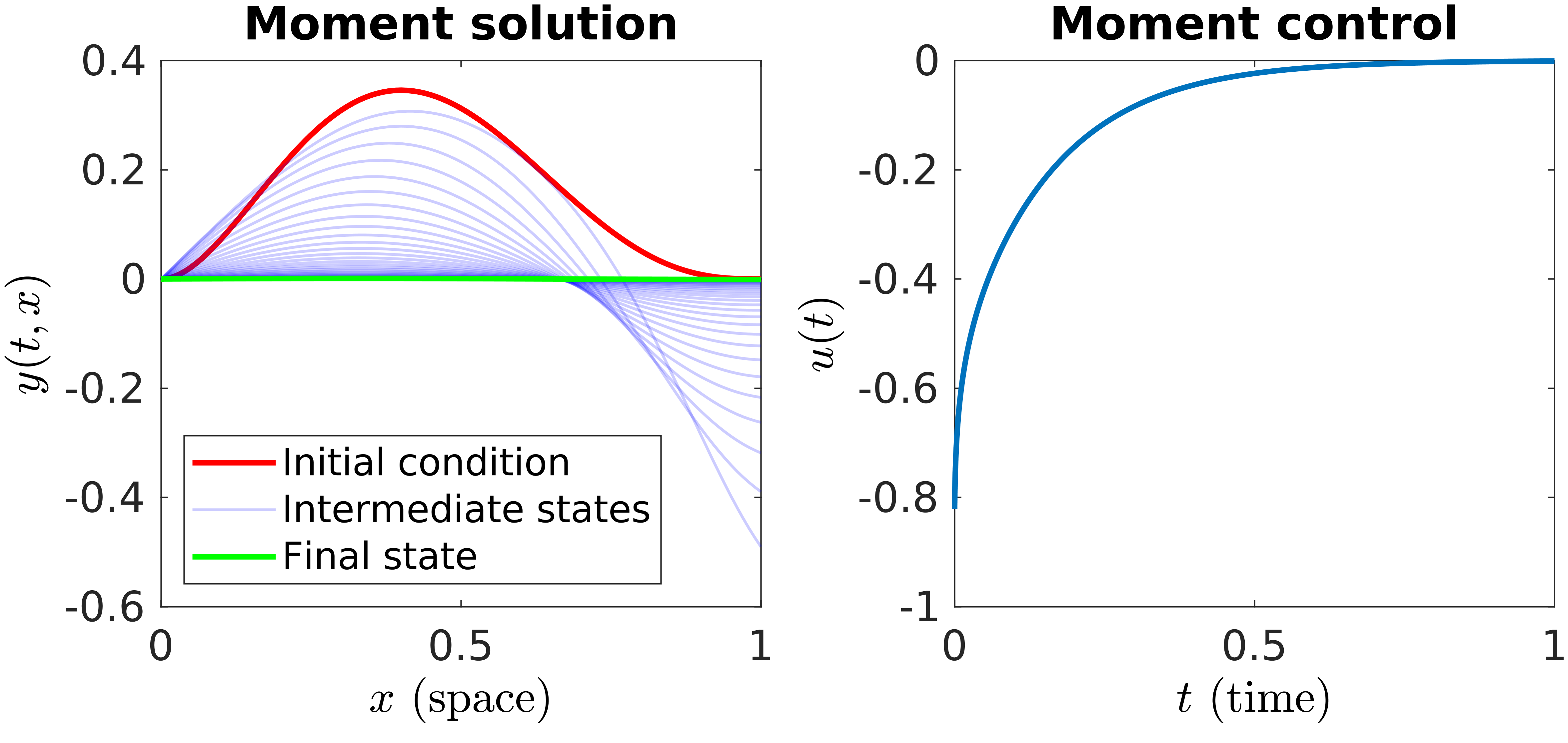}
    \caption{Numerical simulation of \eqref{eq:heat1D} with $\eta = 0$ and a linear control extracted from the moments.}
    \label{fig:numsim_moment_control_linear_heat}
\end{figure}

In this section, we apply the moment-based method to our problem. The LMI problem is solved at the relaxation degree $d=6$ and we extract a linear controller of the form \eqref{eq:control_linear_integral}, where $\gamma_{\,\ell}$ is chosen to be constant, (i.e., $m=0$) and $p=d-1=5$. The least-squares solution of the linear system \eqref{eq:rectangular_sys} in this case is given by the constant polynomial $\gamma_{\,\ell}(t,x) \approx -4.927$. We present the results of the numerical simulations in Fig. \ref{fig:numsim_moment_control_linear_heat}. We observe that the linear control derived from the pseudo-moments successfully controls the solution to 0.

We can compare the value of the cost function for the pseudo-moment controller $L_{\mathrm{mom}} \approx \num{1.850e-3}$ to the optimal value derived from the LQR controller. It corresponds to a relative error of $1.148$\% from the optimal value $L_{\mathrm{LQR}}$. Thus, the linear controller extracted from the pseudo-moments is close to optimality and could be a promising alternative to the LQR. Indeed, in cases where the number of degrees of freedom $N$ is too large to solve the Riccati equation, the moment-based method comes in handy as it does not rely on spatio-temporal gridding and can provide linear controllers similar to the LQR. This last remark is justified since the LQR can be expressed as a linear integral transform of the solution, where the kernel can be deduced from the adjoint state (see for instance \cite[Chapter 3]{PMP-ref}). If the kernel is regular enough, \eqref{eq:control_linear_integral} is a good approximation of the LQR.

\subsection{The nonlinear case}
\label{sec:nonlinear_case}

\begin{figure}[ht]
    \centering
    \includegraphics[width=0.97\columnwidth]{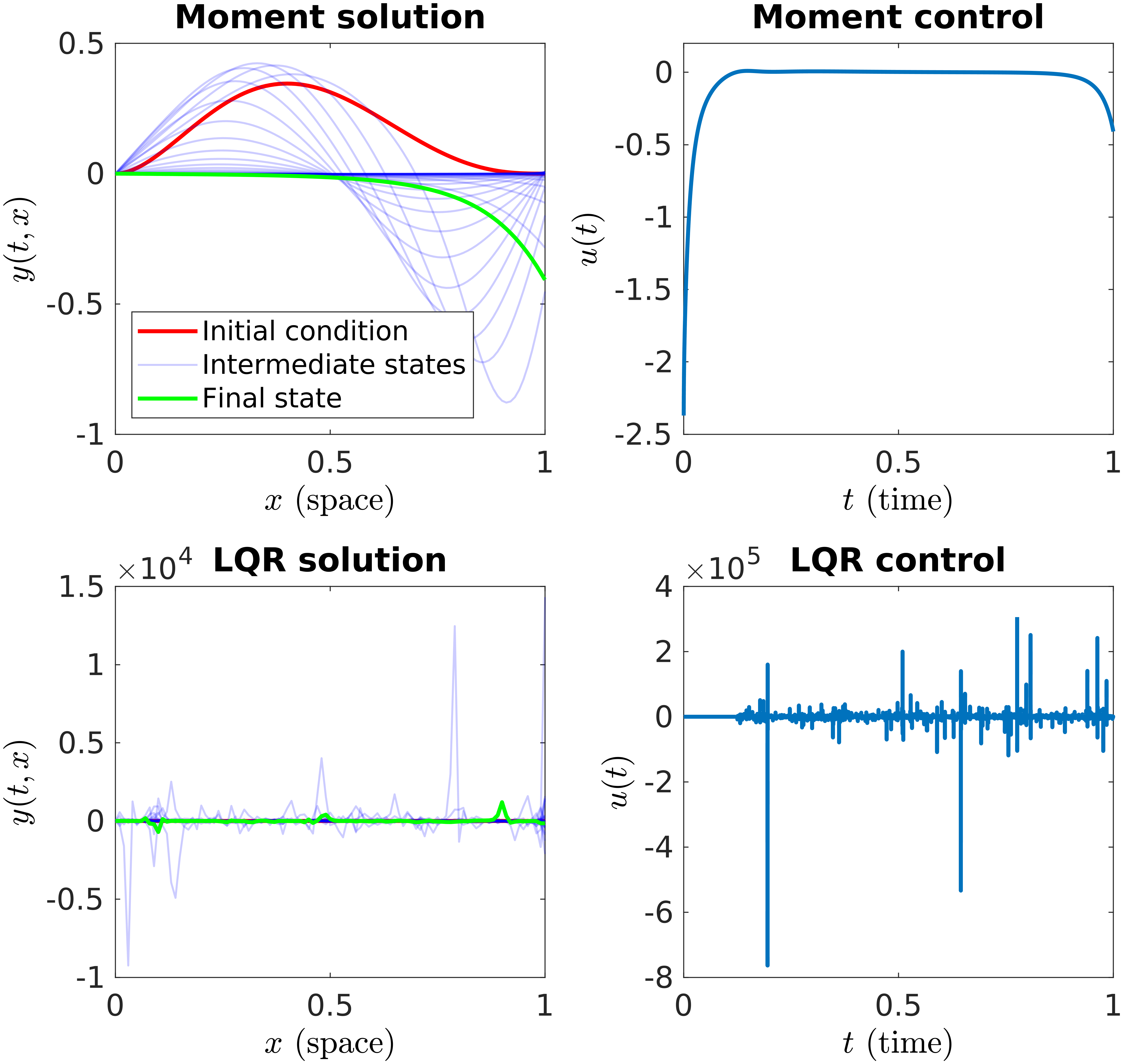}
    \caption{Numerical simulation of \eqref{eq:heat1D} with $\eta \neq 0$. Comparison of a nonlinear control law computed from the moments (top row) with the LQR controller (bottom row).}
    \label{fig:numsim_control_nonlinear_heat}
\end{figure}

We now consider the semilinear PDE where $\eta \neq 0$. The purpose of this section is to demonstrate a scenario where the LQR controller, derived from the linearized equation with $\eta = 0$ in Section \ref{sec:LQR}, fails to control the nonlinear PDE \eqref{eq:heat1D} towards the target state $y=0$. We present the results of the numerical simulations in Fig. \ref{fig:numsim_control_nonlinear_heat}. The parameters used for the simulation are identical to those used for Fig. \ref{fig:numsim_LQR_linear_heat}, with $\eta$ set to $\eta = 13 \alpha$. For this specific set of parameters, numerical simulations of the PDE with null boundary control result in a finite time blow-up of the solution at $t \approx 0.1$. We solve the LMI relaxations for the relaxation degree $d=6$ and compute the control from the pseudo-moments as a nonlinear feedback on the solution of the form \eqref{eq:control_semilinear_integral}. We choose this specific form because it usually provides better performance results compared to the most general form \eqref{eq:control_nonlinear_integral}. The combination of parameters that produced the best results is $r = 3$, with $\gamma_{s\ell}$ being a first-degree polynomial and $\delta_{s\ell}$ being a constant.
The resolution of the linear system \eqref{eq:rectangular_sys} yields
\begin{subequations}
    \begin{align}
        \gamma_{\,s\ell}(t,x) &\approx -15.005 + 21.374 \, t + 1.231 \, x, \label{eq:gamma_expression}\\
        \delta_{\,s\ell}(t,x) &\approx 3.369. \label{eq:delta_expression}
    \end{align}
\end{subequations}
%\CL{We observe that, even if the LQR controller still manages to delay the blow-up outside of the time interval $[0,1]$, it fails to steer the solution towards 0. On the contrary, the semilinear controller derived from the moments successfully stabilizes the solution, at least for $t \leq 0.9$. When $t$ approaches 1, the moment solution starts moving away from the target state 0. Despite this behavior, we consider that this semilinear controller is satisfying enough as it drives the solution sufficiently close to 0 for standard linear control techniques to take over and maintain the solution near the target.}
%\EF{test}
We observe that the controller derived from the moments successfully steers the solution towards 0 at least on the time window $[0,0.9]$, whereas the LQR controller does not.  For $t\geq 0.9$, the moment controller starts driving the solution away from 0.  However, since a standard linear controller could take over when the solution is close to 0, the behavior of the moment controller is satisfactory. These results are promising because the method (cf. Fig. \ref{fig:steps_extract_controller}) can be readily applied to semilinear heat equations with polynomial nonlinear terms other than cubic. A trial and error type of search on the parameters $m$ and $p$ can be easily implemented to extract the nonlinear controller that gives the best results.

\addtolength{\textheight}{-0.0486cm}

\section{CONCLUSIONS}
\label{sec:conclusion}
In this paper, an extended moment-SOS formulation, as well as an associated controller design method for a boundary optimal control problem of a 1D semilinear heat equation are presented. The proposed framework is an extension of that of \cite{KHL-2018} and allows one to consider a quadratic cost on the control. Additionally, a nonlinear feedback controller acting on the boundary of the domain is derived from pseudo-moments as the integral over the whole domain of a multivariate polynomial. The efficiency of our approach is demonstrated through numerical experiments in a case where traditional methods based on the linearized PDE fail because of the nonlinearity of the problem.

Future work will pursue several research directions including: (a) changing the polynomial basis in the moment-SOS hierarchy because the monomial basis is usually badly conditioned (e.g., the Chebyshev polynomial basis as in \cite{chebyshev-ref}); (b) using the Christoffel-Darboux kernel to extract controllers from pseudo-moments (see for instance \cite{marx-darboux} for the reconstruction of the solution of a nonlinear PDE); (c) extending the method to 2D semilinear heat equations and, as a long term goal, to fluid-structure interaction problems.

\bibliography{root.bib}%biblio/references}

\begin{thebibliography}{10}
\providecommand{\url}[1]{#1}
\csname url@rmstyle\endcsname
\providecommand{\newblock}{\relax}
\providecommand{\bibinfo}[2]{#2}
\providecommand\BIBentrySTDinterwordspacing{\spaceskip=0pt\relax}
\providecommand\BIBentryALTinterwordstretchfactor{4}
\providecommand\BIBentryALTinterwordspacing{\spaceskip=\fontdimen2\font plus
\BIBentryALTinterwordstretchfactor\fontdimen3\font minus \fontdimen4\font\relax}
\providecommand\BIBforeignlanguage[2]{{%
\expandafter\ifx\csname l@#1\endcsname\relax
\typeout{** WARNING: IEEEtran.bst: No hyphenation pattern has been}%
\typeout{** loaded for the language `#1'. Using the pattern for}%
\typeout{** the default language instead.}%
\else
\language=\csname l@#1\endcsname
\fi
#2}}

\bibitem{control-navier-stokes}
E.~Bänsch, P.~Benner, J.~Saak, and H.~K. Weichelt, ``\BIBforeignlanguage{en}{Riccati-based boundary feedback stabilization of incompressible {Navier}--{Stokes} flows},'' \emph{\BIBforeignlanguage{en}{SIAM Journal on Scientific Computing}}, vol.~37, no.~2, pp. A832--A858, 2015.

\bibitem{fsi-1}
A.~Mini, C.~Lerch, R.~Wüchner, and K.-U. Bletzinger, ``Computational closed-loop control of fluid-structure interaction ({FSCI}) for lightweight structures,'' \emph{PAMM}, vol.~16, no.~1, pp. 15--18, 2016.

\bibitem{num-method-fsi}
M.~Fourni\'{e}, M.~Ndiaye, and J.-P. Raymond, ``Feedback stabilization of a two-dimensional fluid-structure interaction system with mixed boundary conditions,'' \emph{SIAM Journal on Control and Optimization}, vol.~57, no.~5, pp. 3322--3359, 2019.

\bibitem{course-LQR}
A.~Bensoussan, G.~Prato, M.~C. Delfour, and S.~K. Mitter, Eds., \emph{\BIBforeignlanguage{en}{Representation and {Control} of {Infinite} {Dimensional} {Systems}}}, 2nd~ed., ser. Systems \& {Control}: {Foundations} \& {Applications}.\hskip 1em plus 0.5em minus 0.4em\relax Boston, MA: Birkhäuser Boston, 2007.

\bibitem{riccati-proj}
C.~Airiau, J.-M. Buchot, R.~K. Dubey, M.~Fournié, J.-P. Raymond, and J.~Weller-Calvo, ``\BIBforeignlanguage{en}{Stabilization and best actuator location for the {Navier}--{Stokes} equations},'' \emph{\BIBforeignlanguage{en}{SIAM Journal on Scientific Computing}}, vol.~39, no.~5, pp. B993--B1020, 2017.

\bibitem{quarteroni-optim-then-discrete}
A.~Quarteroni, \emph{\BIBforeignlanguage{en}{Numerical {Models} for {Differential} {Problems}}}.\hskip 1em plus 0.5em minus 0.4em\relax Milano: Springer Milan, 2014.

\bibitem{backstepping-ref}
M.~Krstic and A.~Smyshlyaev, \emph{Boundary Control of PDEs}.\hskip 1em plus 0.5em minus 0.4em\relax Philadelphia, PA: Society for Industrial and Applied Mathematics, 2008.

\bibitem{coron-CLF}
J.-M. Coron, \emph{\BIBforeignlanguage{en}{Control and {Nonlinearity}}}, ser. Mathematical {Surveys} and {Monographs}.\hskip 1em plus 0.5em minus 0.4em\relax Providence, Rhode Island: American Mathematical Society, 2009, vol. 136.

\bibitem{DMD-ref}
H.~Arbabi, M.~Korda, and I.~Mezić, ``A data-driven {Koopman} model predictive control framework for nonlinear partial differential equations,'' in \emph{2018 IEEE Conference on Decision and Control (CDC)}, 2018, pp. 6409--6414.

\bibitem{peitz2019koopman}
S.~Peitz and S.~Klus, ``Koopman operator-based model reduction for switched-system control of {PDEs},'' \emph{Automatica}, vol. 106, pp. 184--191, 2019.

\bibitem{PMP-ref}
F.~Tröltzsch, \emph{\BIBforeignlanguage{en}{Optimal Control of Partial Differential Equations: Theory, Methods, and Applications}}, ser. Graduate studies in mathematics.\hskip 1em plus 0.5em minus 0.4em\relax Providence, R.I: American Mathematical Society, 2010, no. v. 112.

\bibitem{HJB-2}
D.~Kalise, S.~Kundu, and K.~Kunisch, ``Robust feedback control of nonlinear {PDEs} by numerical approximation of high-dimensional {Hamilton}--{Jacobi}--{Isaacs} equations,'' \emph{SIAM Journal on Applied Dynamical Systems}, vol.~19, no.~2, pp. 1496--1524, 2020.

\bibitem{SDRE-ref}
A.~Alla, D.~Kalise, and V.~Simoncini, ``\BIBforeignlanguage{en}{State-dependent {Riccati} equation feedback stabilization for nonlinear {PDEs}},'' \emph{\BIBforeignlanguage{en}{Advances in Computational Mathematics}}, vol.~49, no.~1, p.~9, 2023.

\bibitem{machine-learning}
T.~Duriez, S.~L. Brunton, and B.~R. Noack, \emph{\BIBforeignlanguage{en}{Machine {Learning} {Control} – {Taming} {Nonlinear} {Dynamics} and {Turbulence}}}, ser. Fluid {Mechanics} and {Its} {Applications}.\hskip 1em plus 0.5em minus 0.4em\relax Springer International Publishing, 2017, vol. 116.

\bibitem{KHL-2018}
M.~Korda, D.~Henrion, and J.~B. Lasserre, ``Moments and convex optimization for analysis and control of nonlinear {PDEs},'' in \emph{Numerical Control: Part A}, ser. Handbook of Numerical Analysis, E.~Trélat and E.~Zuazua, Eds.\hskip 1em plus 0.5em minus 0.4em\relax Elsevier, 2022, vol.~23, pp. 339--366.

\bibitem{moment-entropy-PDEs}
S.~Marx, T.~Weisser, D.~Henrion, and J.~B. Lasserre, ``A moment approach for entropy solutions to nonlinear hyperbolic {PDEs},'' \emph{MCRF}, vol.~10, no.~1, pp. 113--140, 2020.

\bibitem{moment-entropy-param-dependent}
C.~Cardoen, S.~Marx, A.~Nouy, and N.~Seguin, ``A moment approach for entropy solutions of parameter-dependent hyperbolic conservation laws,'' \emph{Numerische Mathematik}, vol. 156, pp. 1--36, 2024.

\bibitem{moment-integral-variational-1}
G.~Valmorbida, M.~Ahmadi, and A.~Papachristodoulou, ``Stability analysis for a class of partial differential equations via semidefinite programming,'' \emph{IEEE Transactions on Automatic Control}, vol.~61, no.~6, pp. 1649--1654, 2016.

\bibitem{moment-integral-variational-2}
A.~Chernyavsky, J.~J. Bramburger, G.~Fantuzzi, and D.~Goluskin, ``Convex relaxations of integral variational problems: Pointwise dual relaxation and sum-of-squares optimization,'' \emph{SIAM Journal on Optimization}, vol.~33, no.~2, pp. 481--512, 2023.

\bibitem{valmorbida_reviewer}
A.~Gahlawat and G.~Valmorbida, ``A semi-definite programming approach to stability analysis of linear partial differential equations,'' in \emph{2017 IEEE 56th Annual Conference on Decision and Control (CDC)}, 2017, pp. 1882--1887.

\bibitem{EDP-magron}
V.~Magron and C.~Prieur, ``Optimal control of {PDEs} using occupation measures and {SDP} relaxations,'' \emph{IMA Journal of Mathematical Control and Information}, 2017.

\bibitem{moment-inf-dim}
D.~{Henrion}, M.~{Infusino}, S.~{Kuhlmann}, and V.~{Vinnikov}, ``{Infinite-dimensional moment-SOS hierarchy for nonlinear partial differential equations},'' \emph{arXiv e-prints}, p. arXiv:2305.18768, 2023.

\bibitem{well_posed_pde}
D.~Pighin and E.~Zuazua, ``Controllability under positivity constraints of semilinear heat equations,'' \emph{Mathematical Control and Related Fields}, vol.~8, no. 3\&4, pp. 935--964, 2018.

\bibitem{korda2022gap}
M.~Korda and R.~Rios-Zertuche, ``The gap between a variational problem and its occupation measure relaxation,'' \emph{arXiv preprint arXiv:2205.14132}, 2022.

\bibitem{henrion2024occupation}
D.~Henrion, M.~Korda, M.~Kruzik, and R.~Rios-Zertuche, ``Occupation measure relaxations in variational problems: the role of convexity,'' \emph{SIAM Journal on Optimization}, vol.~34, no.~2, pp. 1708--1731, 2024.

\bibitem{fantuzzi2022sharpness}
G.~Fantuzzi and I.~Tobasco, ``Sharpness and non-sharpness of occupation measure bounds for integral variational problems,'' \emph{arXiv preprint arXiv:2207.13570}, 2022.

\bibitem{book-lasserre}
J.~B. Lasserre, \emph{Moments, Positive Polynomials and Their Applications}.\hskip 1em plus 0.5em minus 0.4em\relax Imperial College Press, 2009.

\bibitem{glopti-article}
D.~Henrion, J.~B. Lasserre, and J.~Löfberg, ``{GloptiPoly} 3: Moments, optimization and semidefinite programming,'' \emph{Optimization Methods and Software}, vol.~24, 2007.

\bibitem{sra_interior-point_2011}
M.~Andersen, J.~Dahl, Z.~Liu, and L.~Vandenberghe, ``\BIBforeignlanguage{en}{Interior-point methods for large-scale cone programming},'' in \emph{\BIBforeignlanguage{en}{Optimization for {Machine} {Learning}}}, S.~Sra, S.~Nowozin, and S.~J. Wright, Eds.\hskip 1em plus 0.5em minus 0.4em\relax The MIT Press, 2011, pp. 55--84.

\bibitem{ARE-book}
P.~Lancaster and L.~Rodman, \emph{{Algebraic Riccati Equations}}.\hskip 1em plus 0.5em minus 0.4em\relax Oxford University Press, 1995.

\bibitem{chebyshev-ref}
D.~Henrion, ``\BIBforeignlanguage{eng}{Semidefinite characterisation of invariant measures for one-dimensional discrete dynamical systems},'' \emph{\BIBforeignlanguage{eng}{Kybernetika}}, vol.~48, no.~6, pp. 1089--1099, 2012.

\bibitem{marx-darboux}
S.~Marx, E.~Pauwels, T.~Weisser, D.~Henrion, and J.~B. Lasserre, ``Semi-algebraic approximation using {Christoffel--Darboux} kernel,'' \emph{Constructive Approximation}, vol.~54, pp. 391--429, 2021.

\end{thebibliography}

\end{document}